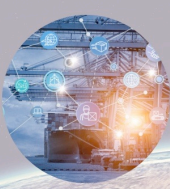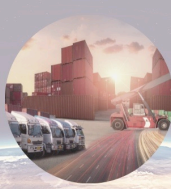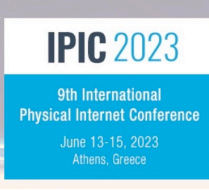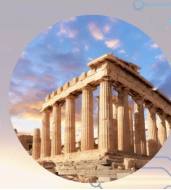

# Stochastic Service Network Design with Different Operational Patterns for Hyperconnected Relay Transportation

Jingze Li, Xiaoyue Liu, Mathieu Dahan, and Benoit Montreuil

Physical Internet Center, H. Milton Stewart School of Industrial & Systems Engineering,
Georgia Institute of Technology, Atlanta, GA 30332
Corresponding author: jonyli@gatech.edu

***Abstract:*** *Hyperconnected relay transportation enables using a relay system of short-haul drivers to deliver long-haul shipments collectively, which helps address root causes of trucker shortage issues by transforming working conditions with potentials of daily returning home, accessing consistent schedules, and facilitating load matching. This paper investigates hyperconnected relay transportation as a sustainable solution to trucker shortage issues through a logistics platform. We propose a two-stage programming model to optimize consistent working schedules for short-haul drivers while minimizing transportation costs. The first stage involves opening services and contracting truckers under demand uncertainty, where each service has a service route and approximate service schedules adhering to USA federal short-haul hour-of-service regulations. The second stage assigns hauling capacities to open services and manages commodity shipping or outsourcing given the demand realization. We extend the model formulation to account for various operational patterns (e.g., freight loading and unloading or hauler swapping) and schedule consistency requirements (e.g., weekly or daily consistency). A scenario-based approach is employed to solve the model for a case study of automotive delivery in the Southeast USA region. The experimental results validate the proposed approach, and further explore the impact of stochastic demands, operational patterns, consistent schedules, and hauling capacities on hyperconnected service network design. This research aims to offer practical guidance to practitioners in the trucking industry.*

***Keywords:*** *Hyperconnected Relay Transportation; Logistics Platform; Stochastic Service Network Design; Short-Haul Truckers; Hour-of-Service Regulations; Demand Uncertainty; Operational Patterns; Consistent Schedules; Hauling Capacities; Physical Internet*

***Conference Topic(s):*** *Interconnected Freight Transport*

***Physical Internet Roadmap*** (Link): Select the most relevant area for your paper: ☐ PI Nodes, ☒ PI Networks, ☐ System of Logistics Networks, ☐ Access and Adoption, ☐ Governance.

## 1 Introduction

Transportation is crucial in moving people and goods across the world, and over the years, various transportation systems have emerged to meet the diverse needs of society. Currently, the prevalent system combines point-to-point and hub-and-spoke transport. The former involves direct movement from one location to another, whereas the latter utilizes a series of spokes to connect hubs that function as transfer points. Hyperconnected relay transport is a recent concept stemming from the Physical Internet (PI) vision (Montreuil, 2011). Inspired by the Digital Internet, where data packets independently traveling through routers and cables form a complete message upon arrival, it relies on a meshed network of interconnected relay nodes for multi-segment and intermodal transport. Unlike hub-and-spoke transport, hyperconnected relay transport enables each hub to act as a transfer point, achieving flexible consolidation opportunities and delivery options. Additionally, it enables a relay system of short-haul drivers

Jingze Li, Xiaoyue Liu, Mathieu Dahan, Benoit Montreuil

to collectively transport long-haul freight, addressing the core issues of trucker shortages by transforming working conditions with benefits such as daily home returns, consistent schedules, and improved driver-freight matching. In this paper, we want to provide a sustainable solution to trucker shortage from the perspective of hyperconnected relay transport.

In recent years, the global logistics industry has witnessed a rise in the number of logistics platforms driven by the growth of e-commerce and the demand for fast, reliable, and cost-effective logistics solutions. These platforms offer numerous advantages, including streamlined market access, efficient load matching, enhanced shipment visibility, and increased delivery efficiency. This paper explores the application of a logistics platform to facilitate the implementation of hyperconnected relay transport. Specifically, it considers a scenario where a logistics platform manages services over a network of relay hubs. The platform receives long-haul loads from shippers and offers short-haul contracts to carriers and owner-operator truckers. The primary goals are to deliver the loads efficiently and reliably within shippers' requested time windows and to improve truckers' working conditions by getting them back home daily. The platform offers compelling value propositions to its three main stakeholders. Firstly, shippers benefit from the ability to express their transportation needs for the foreseeable future, with consistent access to the required transport capacity at lower costs and a wider range of delivery options. Secondly, carriers attain opportunities to secure contracts, revenues, and shipments for their truckers well in advance. Lastly, the platform ensures both company-employed and owner-operator truckers the ability to return home each day through shorter-haul routes, increased visibility into upcoming tasks, and more consistent and predictable schedules.

In this paper, we propose a methodology for optimizing the platform's tactical decisions to persistently achieving its goals and value propositions. We focus on designing a truck-based hyperconnected service network by solving a two-stage stochastic programming model. The objective is to create optimal consistent working schedules for contracted short-haul truckers considering stochastic demands, while minimizing the total transportation cost for the platform. In the first stage, we make decisions regarding which services to open and determine the number of truckers to contract for each open service, taking into account uncertain demands. Each service has a defined route and approximate schedules, complying with USA federal short-haul hour-of-service regulations. In the second stage, after the realization of the demand scenario, the decisions are to determine whether to ship or outsource commodities, how to ship them through open services, and how to assign hauling capacities to each open service. Additionally, we extend the model to incorporate different operational patterns, such as freight loading and unloading or hauler swapping, as well as schedule consistency requirements (e.g., weekly or daily consistency). We solve the model using a scenario-based approach for an automotive delivery case in the Southeast USA region as a testbed. The transportation involves carriers and truckers responsible for delivering vehicles from multiple Original Equipment Manufacturers (OEMs), railheads, and ports to dealers through a relay hub network. The results demonstrate that our proposed approaches can significantly improve the drivers' work-life balance by allowing short-haul drivers to return home daily and maintaining steady working schedules, while considering short-haul hour-of-service regulations and delivery timeliness. We also conduct comparisons of the service network design under various operational patterns, schedule consistency requirements and hauling capacity options to analyze their impact on the network structure and capability.

The full paper is organized as follows. Section 2 discusses the related literature. Section 3 proposes the two-stage stochastic model formulation of hyperconnected service network design and explores its variants. Section 4 analyzes the computational results. Section 5 summarizes the contributions, limitations, and future work directions.

## 2 Related Literature

In the context of hyperconnected relay transportation, researchers have conducted comprehensive assessments of its performance from economic, environmental, and societal perspectives through simulation-based experiments (Hakimi et al., 2012; Sarraj et al., 2014; Hakimi et al., 2015). By employing case studies with large-scale industrial data, they have demonstrated the substantial improvements on efficiency and sustainability by implementing hyperconnected relay transportation. These improvements encompass a range of factors, including reduced $CO_2$ emissions, cost savings, decreased lead time, improved delivery travel time and so forth. In addition to the aforementioned assessment research, which implements and operates hyperconnected relay transport through heuristic protocols, researchers have also delved into solution design research to address key planning and operational decisions induced by the concept of the Physical Internet in hyperconnected relay transport (Pan et al., 2017). Orenstein et al. (2022) developed a mathematical heuristic for routing and scheduling vehicles involved in parcel transfer, as well as a parcel routing mechanism within the hyperconnected service network. Their simulation study demonstrated the effectiveness and advantages of the proposed approach compared to a tree-like service. Li et al. (2022) designed an operating system based on a multi-agent architecture to tackle daily large-scale operational decision-making related to generating shipments, coordinating shipments, tractors, and trailers, as well as assigning and scheduling truckers. Their simulation results outperformed conventional end-to-end transportation in terms of delivery timeliness, driver at-home time, and total operational cost. Qiao et al. (2016) investigated a dynamic pricing model based on an auction mechanism to optimize carrier bid prices. They considered PI-hubs as spot freight markets where less-than-truckload requests arrive over time for short durations. These studies focus on optimizing routing, scheduling, and pricing decisions with the knowledge of demands in hyperconnected relay transport. In this paper, we address a significant difference by considering how to plan logistics services and make contracts with carriers before the knowledge of demands, given the novel business context of a logistics platform. Additionally, we explicitly address factors such as hub operations patterns, schedule consistency, and hauling capacities, which have not been extensively explored in previous literature. We present model variants and experimental results to analyze their impacts on hyperconnected service network design.

Another related literature topic is called service network design, which involves planning routing and scheduling of services and shipments through a network of terminals. Many researchers have approached the modeling of the service network design problem by utilizing the time-space network formulation and incorporating customized rules for various settings, including the network infrastructure, transportation operations, and fleet composition (Scherr et al., 2019; Medina et al., 2019). A recent focus in service network design is addressing uncertainty to enhance robustness and stability, known as stochastic service network design. Two common sources of uncertainty are demands and traffic time, where the former typically arise prior to transportation activities and the latter occurs during and after these activities. Bai et al. (2014) introduced rerouting as a flexible approach for freighters to adapt to demand uncertainty, in addition to outsourcing. Wang et al. (2016) considered variable service capacities and ad-hoc handling to tackle stochastic demands. Lanza et al. (2021) compared travel time uncertainty to demand uncertainty and proposed a model that explicitly incorporates travel time uncertainty and quality targets. In our study, we focus on demand uncertainty in developing consistent approximate schedules, referred to as services, for contracted short-haul truckers. We represent this uncertainty through cyclic scenarios of demand patterns and formulate a two-stage model. Such "inherently two-stage problem" focuses on the first stage of designing a network and creates a correct understanding of how the network will be operated through the second stage. It actually simplifies the multi-stage nature of the real problem such that we can bypass intricate details that are more pertinent to the dynamic operational phase. We also list refining approximate schedules by accounting for stochastic traffic time as one of our future works. Such idea of approximation-then-refining is inspired by Boland et al. (2017), who introduced a systematic computational method known as dynamic discretization discovery for the continuous service network design in a general context.

## 3 Methodology

We focus on the case where a platform manages the logistics service over a provided relay hub network. The physical network is denoted by $\mathcal{G}^P = (\mathcal{N}^P, \mathcal{A}^P)$ with node set $\mathcal{N}^P$ representing hub nodes and arc set $\mathcal{A}^P$ representing connected arcs between hub nodes. A planning horizon is considered and discretized into $T + 1$ evenly distributed time instants, which is denoted as $\mathcal{T} = \{0,1,2,\ldots,T\}$. We then construct a time-space network $\mathcal{G} = (\mathcal{N}, \mathcal{A})$ based on the physical network $\mathcal{G}^P = (\mathcal{N}^P, \mathcal{A}^P)$ and the discretized planning horizon $\mathcal{T}$. The node set $\mathcal{N}$ is attained by replicating each node in $\mathcal{N}^P$ for $T + 1$ times, where $\mathcal{N} = \{(n,t) | n \in \mathcal{N}^P, t \in \mathcal{T}\}$. The arc set $\mathcal{A}$ consists of a moving arc set $\mathcal{A}^M$ and a holding arc set $\mathcal{A}^H$. The former includes the arcs between time replicates of two different hub nodes to model the movements of freight or truckers, while the latter incorporates the arcs between two replicates of the same hub node at two consecutive time instants, used for modeling idle time or processing time of freight or truckers. Each arc $a \in \mathcal{A}$ has the form $a = \big((n_a^1, t_a^1), (n_a^2, t_a^2)\big)$ where $(n_a^1, t_a^1), (n_a^2, t_a^2) \in \mathcal{N}$. Specifically, if $a$ is a moving arc, we will have $n_a^1 \neq n_a^2$ and $t_a^2 - t_a^1$ representing the travel time plus buffer time for the movement from hub node $n_a^1$ to hub node $n_a^2$; if $a$ is a holding arc, we will have $n_a^1 = n_a^2$ and $t_a^2 = t_a^1 + 1$.

The platform makes contracts over a set of potential services $\mathcal{S}$ to transport freight. Each service $s \in \mathcal{S}$ has a fixed service route and approximate service schedules, represented by $s = \{a_s^1, \ldots, a_s^{r_s}\}$, where $a_s^i \in \mathcal{A}^M$ indicates the $i$-th movement of service $s$ for $i \in \{1, \ldots, r_s\}$. For example, a service $s = \{((n_1, t_1), (n_2, t_2)), ((n_2, t_3), (n_1, t_4))\}$ represents a service planned to start from hub node $n_1$ at time instant $t_1$, arrive at hub node $n_2$ at time instant $t_2$, then leave hub node $n_2$ at time instant $t_3$, and return to hub node $n_1$ at time instant $t_4$. The platform needs to decide which service to open and how many truckers to contract over each open service before the knowledge of demands. Once truckers are contracted, the contract fees are paid and cannot be cancelled. Each service $s \in \mathcal{S}$ has a contacted trucker capacity $q_s$. Once the truckers are contracted to services, the platform can then add a hauling capacity $u \in \mathcal{U}$ (or zero capacity in case of excess capacity) to each contracted trucker based on the actual demand scenario.

The platform receives the transportation requests for multiple commodities. Each commodity $k \in \mathcal{K}$ has an origin hub $o_k$, a destination hub $d_k$, an entry time $t_k^e$, a due time $t_k^d$, and volume $v_k$. For each commodity, the platform has the option of transporting it by contracted services or outsourcing it to third-party logistics. All commodities are expected to be delivered on time.

We assume that the hyperconnected service network of the logistics platform adheres to a distinct hub operations pattern and maintains certain schedule consistency. Two operational patterns are taken into consideration for hub operations: freight loading and unloading (FLU) and hauler swapping (HS). For FLU, drivers stay with their trucks (including tractors and haulers), and freight can be loaded and unloaded at each hub for crossdocking. On the other hand, for HS, haulers can be separated from drivers and tractors. Once freight is loaded into a hauler, it remains inside the hauler until reaching its destination, without any additional loading or unloading during the delivery process. Furthermore, for FLU, we consider a condition whether freight of the same commodity will travel along a unique commodity path from origin to destination. Such consideration is due to the facts that customers may prefer to receive the products as a whole pack and keep better track of the shipments. Consequently, we have three operational patterns to consider: *(i-1)* freight loading and unloading with multiple commodity paths (FLU - MCP), *(i-2)* freight loading and unloading with single commodity path (FLU - SCP), and hauler swapping (HS). Schedule consistency is also taken into account due to drivers' preference for consistent schedules, such as weekly or daily schedules. We denote the planning horizon as $\mathcal{T}$, which is composed of multiple schedule cycles, i.e., $\mathcal{T} = \mathcal{T}^1 \cup ... \cup \mathcal{T}^C$, where $C$ represents the total number of cycles. Each service $s$ has its start time determined by the service cycle $c_s$ and the specific start time $t_s^c$ within that cycle.

## 3.1 Two-Stage Stochastic Programming Formulation

We develop a two-stage stochastic programming model for the platform to design the hyperconnected logistics service network. We first provide the formulation with operational pattern as FLU - MCP and schedule consistency with total number of cycles $C = 1$. For FLU - MCP, a trucker refers to a driver, to whom we can assign trucks with different hauling capacities. The objective is to minimize total transportation cost consisting of trucker contract cost and hauler rental cost, while building optimal consistent service schedules for short-haul relay truckers. The first-stage decision variable is $X_s$, an integer variable representing the number of drivers to contract to short-haul service $s \in \mathcal{S}$. The second-stage variables consist of $Z_k(w)$, $F_{ka}(w)$ and $Y_{su}(w)$. $Z_k(w)$ is a binary variable to indicate whether to outsource commodity $k$ in demand scenario $w$. $F_{ka}(w)$ is a nonnegative continuous variable meaning the volume of commodity $k$ that will traverse arc $a$ in demand scenario $w$. $Y_{su}(w)$ is a nonnegative integer variable representing the number of trucks with capacity size $u$ that will be assigned to service $s$ in demand scenario $w$.

The mathematical formulation is shown through equation (1) – (5). The objective function (1) is to minimize the sum of total contracted cost of drivers as well as the total expected costs of truck rentals and commodity outsourcing with regards to demand uncertainty. Constraint (2) ensures the contracted number of drivers does not exceed the maximal contracted trucker capacity for each service. Constraint (3) assigns trucks of different capacity sizes to contracted drivers. Constraint (4) guarantees the total truck volume capacity is no smaller than the total commodity volume on each arc. Constraint (5) is the commodity flow balance constraint, which also ensures the timely delivery.

$$\min \sum_{s \in \mathcal{S}} c_s^f X_s + E_{w \in \mathcal{W}}\left[\sum_{s \in \mathcal{S}, u \in \mathcal{U}} c_{su}^v Y_{su}(w) + \sum_{k \in \mathcal{K}} c_k^o Z_k(w)\right] \quad (1)$$

$s.t.$

$$0 \leq X_s \leq q_s \quad \forall s \in \mathcal{S} \quad (2)$$

$$\sum_{u \in \mathcal{U}} Y_{su}(w) \leq X_s \quad \forall s \in \mathcal{S}, w \in \mathcal{W} \quad (3)$$

$$\sum_{s \in \mathcal{S}_a, u \in \mathcal{U}} u Y_{su}(w) \geq \sum_{k \in \mathcal{K}} F_{ka}(w) \quad \forall a \in \mathcal{A}^M, w \in \mathcal{W} \quad (4)$$

$$\sum_{a \in A: n_a^2 = n} F_{ka}(w) - \sum_{a \in A: n_a^1 = n} F_{ka}(w) = \begin{cases} v_k(w)(Z_k(w) - 1), & if \ n = (o_k, t_k^e) \\ v_k(w)(1 - Z_k(w)), & if \ n = (d_k, t_k^d) \\ 0, otherwise \end{cases}$$

$$\forall k \in \mathcal{K}, n \in \mathcal{N}, w \in \mathcal{W} \quad (5)$$

## 3.2 Model variants

Based on the model formulation in the Section 3.1, we then provide model variants that account for different operational patterns. For FLU – SCP, given the fact that each commodity stays on a unique shipment path, we maintain decision variables $X_s, Y_{su}(w)$ as nonnegative integer variables and $Z_k(w) \in \{0,1\}$ but change $F_{ka}(w)$ to a binary variable meaning whether shipment of commodity $k$ will traverse arc $a$ in demand scenario $w$. The formulation is updated to (1) – (3) plus (4') and (5'), where constraint (4') and (5') adjust the truck capacity constraint and commodity flow balance constraint respectively by considering a unique shipment path for each commodity.

$$\sum_{s \in S_a, u \in U} u Y_{su}(w) \geq \sum_{k \in K} v_k(w) F_{ka}(w) \qquad \forall\, a \in \mathcal{A}^M, w \in \mathcal{W} \quad (4')$$

$$\sum_{a \in A: n_a^2 = n} F_{ka}(w) - \sum_{a \in A: n_a^1 = n} F_{ka}(w) = \begin{cases} Z_k(w) - 1, & \text{if } n = (o_k, t_k^e) \\ 1 - Z_k(w), & \text{if } n = (d_k, t_k^d) \\ 0, & \text{otherwise} \end{cases}$$

$$\forall k \in \mathcal{K}, n \in \mathcal{N}, w \in \mathcal{W} \quad (5')$$

For HS, we refer a trucker to a driver-tractor pair, to which we can assign the haulers with different hauling sizes. The three decision variables $X_s, Y_{ku}(w)$ and $F_{ku}(w)$ become nonnegative variables now. $X_s$ represents number of driver-tractor pairs contracted to service $s$. $Y_{ku}(w)$ represents number of haulers with size $u$ used for commodity $k$ in scenario $w$. Note that one of the subscripts of $Y$ variables are updated from service index $s$ to commodity index $k$, since haulers now stay with freight instead of drivers and tractors. $F_{ku}(w)$ indicates number of haulers with commodity $k$ traversing arc $a$ in scenario $w$. We keep $Z_k(w)$ as the same.

The updated model formulation is given by (1'') – (5''). The objective function (1'') is to minimize the sum of total contracted cost of driver-tractor pairs as well as the total expected cost of hauler rentals and commodity outsourcing with regards to uncertain demands. Constraint (2'') ensures the contracted number of driver-tractor pairs does not exceed the maximal contracted trucker capacity for each service. Constraint (3'') assigns haulers of different capacity sizes to contracted truckers. Constraint (4'') guarantees each commodity is accommodated into haulers with enough volume capacity. Constraint (5'') is the commodity flow balance constraint, which also ensures the timely delivery.

$$\min \sum_{s \in S} c_s^f X_s + E_{w \in \mathcal{W}}\left[\sum_{k \in K, u \in U} c_{ku}^v Y_{ku}(w) + \sum_{k \in K} c_k^o Z_k(w)\right] \quad (1'')$$

$s.t.$

$$0 \leq X_s \leq x_s^{max} \qquad \forall\, s \in \mathcal{S} \quad (2'')$$

$$\sum_{k \in K} F_{ka}(w) \leq \sum_{s \in S_a} X_s \qquad \forall\, a \in \mathcal{A}^M, w \in \mathcal{W} \quad (3'')$$

$$\sum_{u \in U} u Y_{ku}(w) \geq v_k(w)(1 - Z_k(w)) \qquad \forall\, s \in \mathcal{S}, w \in \mathcal{W} \quad (4'')$$

$$\sum_{a \in A: n_a^2 = n} F_{ka}(w) - \sum_{a \in A: n_a^1 = n} F_{ka}(w) = \begin{cases} -\sum_{u \in \mathcal{K}} Y_{ku}(w), & \text{if } n = (o_k, t_k^e) \\ \sum_{u \in \mathcal{K}} Y_{ku}(w), & \text{if } n = (d_k, t_k^d) \\ 0, & \text{otherwise} \end{cases}$$

$$\forall k \in \mathcal{K}, n \in \mathcal{N}, w \in \mathcal{W} \quad (5'')$$

We also consider model variants accounting for different schedule consistency requirements. A strong version of consistency constraint is to have $X_s = X_{s'}$ if service $s$ and $s'$ have the identical route path and cycle time but just belong to different cycles within the planning horizon. A soft version is to put a penalty on the cycle inconsistency and add it into the objective function, which can be measured by the sum of differences in contracted number of service truckers across cycles. In this paper, we consider the strong version. Given the weekly and daily rates of driver payment and fleet rental are different in reality, we aim to examine the impact of consistency on the network structure and cost calculation.

## 4 Results and Discussion

In this section, we present the computational results using the data derived from a real-life automotive delivery case in the Southeastern USA. A hyperconnected relay hub network is used for transportation with 19 hubs and 95 arcs, as illustrated in Figure 1. Most of hubs are interconnected to multiple adjacent hubs, and each arc is designed to accommodate traffic uncertainty with a robust arc travel time of approximate 5.5 hours. This ensures that every short-haul trucker can complete a round trip on each arc and return to domicile within the daily maximal driving window of 11hours, as mandated by the USA federal hour-of-service regulations. The demand sample considered in this analysis spans from January 1, 2020, to January 14, 2020, and consists of a total of 13,903 ordered vehicles. Notably, 94.9% of commodity flow is observed to move from west to east, while the remaining 5.1% moves from east to west, as depicted in Figure 1. This indicates an inherent imbalance in the commodity flows between the two directions, with the hyperconnected service network design predominately catering to the higher volume of commodities moving towards the east.

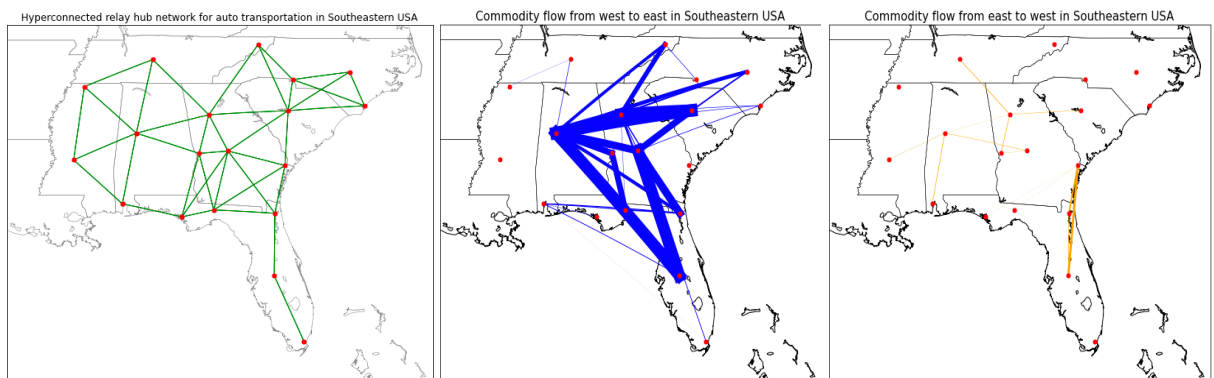

*Figure 1: Hyperconnected hub network, commodity flows towards east (in blue) and west (in orange)*

We utilize a scenario-based approach to solve the hyperconnected service network design model. In our experiments, the planning horizon spans five days, discretized into six-hour time intervals. We generate a total of 30 demand scenarios, encompassing 196 potential commodities between every pair of hubs for the initial four planning days. The delivery time window for these commodities is set as two days. We consider all possible short-haul services with service routes between any two adjacent hubs, while ensuring that service approximate on-duty durations do not exceed the daily maximum on-duty duration of 14 hours imposed by the USA federal hour-of-service regulations. The other parameters for experiments are provided in Table 1. Notably, we employ a consistency cost discount factor of 0.8, which reflects the cost advantages associated with maintaining consistent daily schedules as opposed to weekly consistency.

*Table 1: The other key parameters used for the experiments*

| *Hourly driver contract fee ($)* | 29 | *Hourly size-8 hauler rental fee ($)* | 10 |
|---|---|---|---|
| *Hourly tractor rental fee ($)* | 18 | *Hourly size-4 hauler rental fee ($)* | 5 |
| *Outsourcing cost per vehicle per mile ($)* | 0.93 | *Average mile per hour* | 50 |
| *Contracted trucker capacity per service* | 10 | *Consistency cost discount factor* | 0.8 |

We first conduct two model runs to calculate the Value of the Stochastic Solution (VSS), which is used to measure the importance of incorporating stochasticity. In the first model run, we employ a stochastic model with demand scenarios. In contrast, the second model run involves fixing the first-stage decisions as the deterministic solutions obtained from the deterministic version of the proposed model with demand input as the sample average of demand scenarios,

Jingze Li, Xiaoyue Liu, Mathieu Dahan, Benoit Montreuil

and tests the fixed first-stage decisions in the second stage using demand scenarios. The second model run commonly called "deterministic design in the stochastic model" (Xin et al., 2016). The comparison results of two model runs are shown in Table 2. The deterministic design has less contracted hours of drivers and fewer average rental hours of tractor-hauler pairs than the stochastic one, resulting in 10.3% more average outsourcing rate of commodities. In addition, the VSS is 556,464 – 422,985 = 133,509 in dollars, which shows that the stochastic design can save about 24% of the total expected transportation cost, compared to the deterministic one.

*Table 2: Comparison results of deterministic design vs. stochastic design in a stochastic environment*

| **KPIs \ Demand patterns** | *Deterministic* | *Stochastic* |
|---|---|---|
| *Total contracted hours of drivers (hrs)* | 9,408 | 12,444 |
| *Average rental hours of tractors (hrs)* | 8,023 | 9,285 |
| *Average rental hours of haulers (hrs)* | 8,023 | 9,285 |
| *Average outsourcing rate of commodities* | 10.3% | 0% |
| *Total expected transportation cost ($)* | 556,494 | 422,985 |

To explore the impact of different operational patterns on the hyperconnected service network design, we then perform three model runs with operational patterns as FLU – MCP, FLU – SCP, and HS respectively. In all three model runs, we with consistency requirement set as weekly and fixed hauler size as 8 vehicles per trucker. The comparison results are summarized in Table 3. From FLU-MCP to FLU-SCP, the contracted hours of drivers and rental hours of tractor-hauler pairs increase by 3.4% and 2.9% respectively. All commodities are shipped through the platform logistics, with a 2.1% rise in total expected transportation cost. The reason is that compared with FLU-MCP, FLU-SCP restricts the delivery of each commodity along a single path to ensure freight integrity. Consequently, consolidation opportunities decrease, resulting in a higher number of trucks operating at increased costs. From FLU-SCP to HS, we can observe a 2.6% decrease in total contracted driver hours. This decrease can be attributed to the fact that drivers now remain with their tractors, resulting in higher costs for opening services and contracting truckers. Additionally, the implementation of HS also led to a discernible increase in the outsourcing rate of commodities, accompanied by an overall elevation in the total expected cost. This outcome arises from the absence of crossdocking in the HS setting, necessitating the outsourcing of several commodities as a favorable alternative. However, it is essential to note that HS offers enhanced freight protection by securely maintaining the goods inside the haulers throughout the entire delivery process. Furthermore, it reduces operational efforts by facilitating the efficient swapping of haulers, as opposed to the time-consuming process of loading and unloading freight. While not explicitly modeled in this paper, these aspects present valuable avenues for future research.

*Table 3: Comparison results of stochastic solutions with three different operational patterns*

| **KPIs \ Operational patterns** | *FLU-MCP* | *FLU-SCP* | *HS* |
|---|---|---|---|
| *Total contracted hours of drivers (hrs)* | 12,444 | 12,864 | 12,528 |
| *Average rental hours of tractors (hrs)* | 9,285 | 9,312 | 12,528 |
| *Average rental hours of haulers (hrs)* | 9,285 | 9,312 | 955.2 |
| *Average outsourcing rate of commodities* | 0% | 0% | 1.3% |
| *Total expected transportation cost ($)* | 422,985 | 431,880 | 492,742 |

The subsequent experiment aims to examine the influence of different consistency requirements and hauling capacities on the hyperconnected service network design, with the operational pattern set as FLU – MCP. Two consistency requirements are considered: weekly consistency and daily consistency. The corresponding results are presented in Table 4. Comparing weekly consistency to daily consistency, we observe that daily consistency results in increased contracted hours for drivers and rental hours for trucks, while yields a lower total expected cost.

This outcome can be attributed to the fact that daily-consistent contracts have, on average, lower driver contracted fees and truck rental costs per hour compared to weekly-consistent contracts. Figure 2 visualizes the open services of a weekly-consistent vs. daily-consistent designs. An interesting observation is that, between hub 0 and hub 1, weekly-consistent design has open services on both routes of hub 0 – hub 1- hub 0 and hub 1 – hub 0 – hub 1 on certain days, while daily-consistent design only has open services on the route of hub 1 – hub 0 - hub 1 repeating every day. Furthermore, based on the results in Table 4 and Figure 3, we can observe that various hauling capacities leads to more contracted hours of drivers and rental hours of tractor-hauler pairs, yet contributes to less total expected transportation cost. The reason is that various hauling capacities offer better flexibilities in fleet selection, particularly during periods of low demands, thus the commodities can be shipped in a more cost-effective manner. Additionally, when considering the total expected cost, we find that daily consistency contributes more to cost savings compared to the various hauling capacities. This is because driver contracted fees carry a higher weight in the total expected cost than truck (tractor/hauler) rental costs.

*Table 4: Comparison results of weekly vs. daily consistency requirements in stochastic environment*

| **Consistent patterns** | *Weekly* | | *Daily* | |
|---|---|---|---|---|
| **KPIs \ Hauling capacity** | *Fixed* | *Various* | *Fixed* | *Various* |
| *Total contracted hours of drivers (hrs)* | 12,444 | 12,348 | 13,500 | 13,620 |
| *Average rental hours of tractors (hrs)* | 9,286 | 10,562 | 9,456 | 11,045 |
| *Average rental hours of haulers (hrs)* | 9,286 | 10,562 | 9,456 | 110,45 |
| *Average outsourcing rate of commodities* | 0% | 0% | 0% | 0% |
| *Total expected transportation cost ($)* | 422,985 | 418,253 | 357,840 | 355,152 |

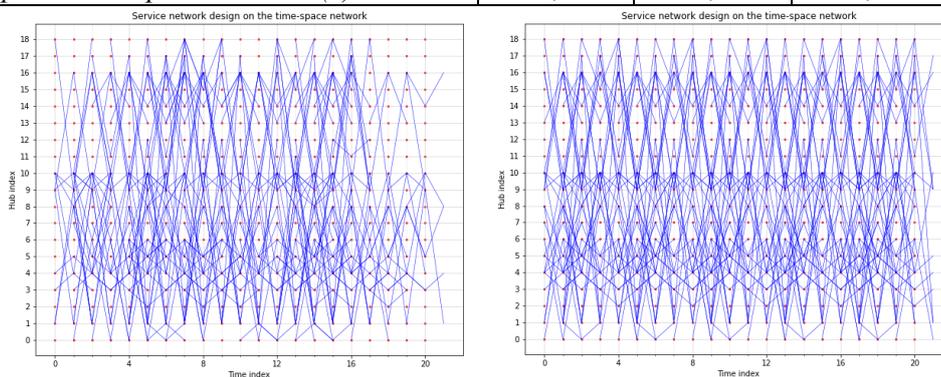

*Figure 2: open services of model designs with weekly consistency (left) vs. daily consistency (right)*

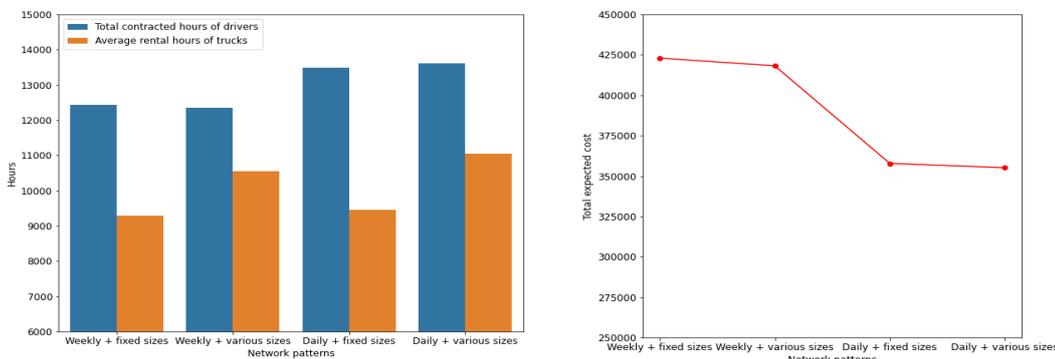

*Figure 3: Total driver contracted hours and truck rental hours (left) and total expected cost (right)*

## 5  Conclusion

The contributions of this paper are threefold. First, it proposes applying hyperconnected relay transportation as a sustainable solution to truck driver shortage issues through a logistics platform as a novel business context. Second, it provides a two-stage stochastic model for hyperconnected service network design of the platform. The formulation supports different operational patterns and schedule consistency requirements. Third, the paper runs the experimental results on an automotive delivery test case in Southeastern USA and discusses the impacts of demand uncertainty, operational patterns, consistent schedules, and various hauling capacities on the service network design, providing guidance to practitioners in industry.

There are several avenues for the future work. The first direction is to develop more advanced computational methods such as bender decomposition or sample average approximation for optimizing larger scale instances. The second is to perform sensitivity analysis upon experimental parameters such as delivery time window and maximal driving time window. The third is to model more route patterns for both short-haul and long-haul, contracted services tailored to trucker preferences, and on-market carrier capacity. The fourth is to refine the approximate service schedules accounting for traffic time stochasticity.

## References

- Montreuil, B. (2011). Toward a Physical Internet: meeting the global logistics sustainability grand challenge. Logistics Research, 3, 71-87.
- Pan, S., Ballot, E., Huang, G. Q., & Montreuil, B. (2017). Physical Internet and interconnected logistics services: research and applications. International Journal of Production Research, 55(9).
- Hakimi, D., Montreuil, B., Sarraj, R., Ballot, E., & Pan, S. (2012, June). Simulating a physical internet enabled mobility web: the case of mass distribution in France. In 9th International Conference on Modeling, Optimization & SIMulation-MOSIM'12 (pp. 10-p).
- Sarraj, R., Ballot, E., Pan, S., Hakimi, D., & Montreuil, B. (2014). Interconnected logistic networks and protocols: simulation-based efficiency assessment. International Journal of Production Research, 52(11), 3185-3208.
- Hakimi, D., Montreuil, B., & Hajji, A. (2015). Simulating Physical Internet Enabled Hyperconnected Semi-Trailer Transportation Systems. In Proceedings of 2nd International Physical Internet Conference, Paris, France.
- Orenstein, I., & Raviv, T. (2022). Parcel delivery using the hyperconnected service network. Transportation Research Part E: Logistics and Transportation Review, 161, 102716.
- Li, J., Montreuil, B., & Campos, M. (2022). Trucker-sensitive Hyperconnected Large-Freight Transportation: An Operating System. In Proceedings of IISE Annual Conference.
- Qiao, B., Pan, S., & Ballot, E. (2019). Dynamic pricing model for less-than-truckload carriers in the Physical Internet. Journal of Intelligent Manufacturing, 30, 2631-2643.
- Scherr, Y. O., Saavedra, B. A. N., Hewitt, M., & Mattfeld, D. C. (2019). Service network design with mixed autonomous fleets. Transportation Research Part E: Logistics and Transportation Review, 124, 40-55.
- Medina, J., Hewitt, M., Lehuédé, F., & Péton, O. (2019). Integrating long-haul and local transportation planning: The service network design and routing problem. EURO Journal on Transportation and Logistics, 8(2), 119-145.
- Bai, R., Wallace, S. W., Li, J., & Chong, A. Y. L. (2014). Stochastic service network design with rerouting. Transportation Research Part B: Methodological, 60, 50-65.
- Wang, X., Crainic, T. G., & Wallace, S. W. (2019). Stochastic network design for planning scheduled transportation services: The value of deterministic solutions. INFORMS Journal on Computing, 31(1), 153-170.
- Lanza, G., Crainic, T. G., Rei, W., & Ricciardi, N. (2021). Scheduled service network design with quality targets and stochastic travel times. European Journal of Operational Research, 288, 30-46.
- Boland, N., Hewitt, M., Marshall, L., & Savelsbergh, M. (2017). The continuous-time service network design problem. Operations research, 65(5), 1303-1321.